\documentclass[11pt]
{amsart}
\usepackage{hyperref}
\hypersetup{nesting=true,debug=true,naturalnames=true}
\usepackage{graphicx,amssymb,upref}

\newtheorem{theorem}{Theorem}

\begin{document}
%
\title{Random walk with barriers on a graph}
%
%
\author{Theo van Uem}
%
\begin{abstract}
 We obtain expected number of arrivals, absorption probabilities and expected time
until absorption for an asymmetric discrete random walk on a graph in
the presence of multiple function barriers. On each edge of the graph and in each vertex (barrier)
specific probabilities are defined. 
\end{abstract}
%
%
%
\maketitle
%
%

\section{Introduction}

Random walk can be used in various disciplines: in economics to model share prices and their
derivatives, in medicine and biology where absorbing barriers give a natural model for a wide variety
of phenomena, in physics as a simplified model of Brownian motion, in ecology to describe individual
animal movements and population dynamics, in statistics to analyze sequential test procedures, in
computer science to estimate the size of the World Wide Web using randomized algorithms.
Burioni and Cassy (2005) give a review of random walks on graphs, where the generalization of the
concept of dimension to inhomogeneous structures, using infinite graphs, is considered.
Durhuus, Jonsson, and Wheater (2006) develop techniques to obtain rigorous bounds on the behavior
of random walks on combs. Using these bounds they calculate the spectral dimension of random
combs with infinite teeth at random positions or teeth with random but finite length.
Random walks have been studied for decades on regular structures such as lattices. We now give a
brief historical review of the use of barriers in a one-dimensional discrete random walk. 
Weesakul (1961)  discussed the classical problem of random walk restricted between a reflecting and an
absorbing barrier. Using generating functions he obtains explicit expressions for the probability of
absorption. Lehner (1963) studies one-dimensional random walk with a partially reflecting barrier
using combinatorial methods. Gupta (1966) introduces the concept of a multiple function barrier
(MFB): a state that can absorb, reflect, let through or hold for a moment. 
Dua, Khadilkar, and Sen (1976) find the bivariate generating functions of the probabilities of a particle reaching a certain state
under different conditions. Percus (1985) considers an asymmetric random walk, with one or two
boundaries, on a one-dimensional lattice. At the boundaries, the walker is either absorbed or reflected
back to the system. Using generating functions the probability distribution of being at position $m$ after
$n$ steps is obtained, as well as the mean number of steps until absorption. El-Shehawey (2000) obtains absorption probabilities at the boundaries for a random walk between one or two partially
absorbing boundaries as well as the conditional mean for the number of steps before stopping given
the absorption at a specified barrier, using conditional probabilities.
In this paper we obtain expected number of arrivals, absorption probabilities and expected time until
absorption for an asymmetric discrete random walk with multiple function barriers. 
Our graph consists of multiple function barriers  (vertices) and states on the edges between the MFB’s.  
On each edge of the graph a random walk with its own  states and
jumping probabilities is introduced. When the walker reaches a multiple function barrier a random
process is activated according to a set of probabilities, or the particle is absorbed in the barrier.
Each barrier has its own probability parameters. In section 2 we use generating
functions to find the expected number of arrivals to any state, the probability of absorption and the
expected time until absorption. In section 3 we analyze some examples of graphs with multiple
function barriers: a star graph and a cycle graph. 

\section{A graph with multiple function barriers}

\subsection{Description of the random walk}

In a finite graph we have vertices $M[0], M[1], \dots,M[N]$ representing the MFB’s. Between $M[i]$ and
$M[j]$ there is a random walk with a finite number of states $n[i,j]$, which we number $1,2,\dots,n[i,j]$ in
the direction from $M[i]$ to $M[j]$ when $i<j$. 
Example 1. Random walk on an interval with two reflecting/absorbing barriers:
\graphicspath{file:///C:/Users/Theo/Desktop/}
\begin{center}
\includegraphics[scale=0.6]{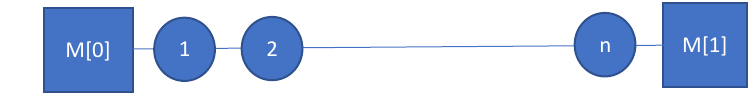} 
\end{center}
Example 2. Random walk on a triangle:
\graphicspath{file:///C:/Users/Theo/Desktop/}
\begin{center}
\includegraphics[scale=0.4]{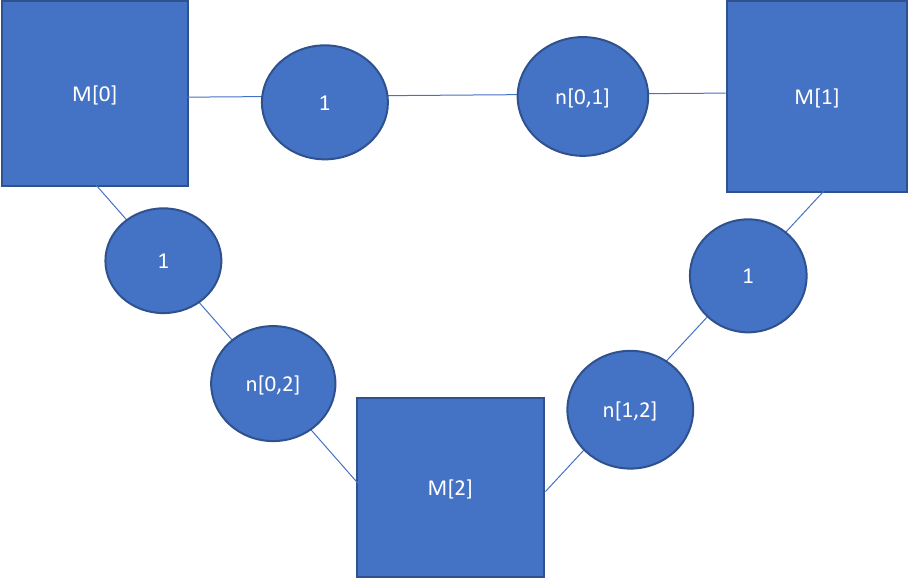} 
\end{center}
We will use the abbreviation $[i,j]$ for the edge between $M[i]$
and $M[j]$. Each random walk from $M[i]$ to $M[j]$ has its own parameters $p=p[i,j]$ and $q=q[i,j]$,
where $p$ is the one-step forward probability and $q$ one-step backward probability $(p+q=1)$. We demand $p[i,j].q[i,j] >0$ for each $i$ and $j$. 
In $M[i]$ there is probability $p_{[i,j]}^{*}$
 to move one step in the direction of $M[j] \ (0\leq i,j\leq N)$ and probability $p_{[i,i]}^{*}$ for absorption in $M[i]$ $(0\leq i \leq N)$, where $\sum_{j=0}^{N}p_{[i,j]}^{*}=1 \ (i=0,1,\dots ,N)$. We start in $M[0]$.
 
 \subsection{Expected number of arrivals}
 We are interested in the expected number of arrivals in the MFB's as well as the expected number of arrivals in the other states of the graph.
 Let $p_{ij}^{(m)}$ be the probability that the system is in state $j$ after $m$ steps when starting in $i$.
 If $j$ is not a MFB:
 \[X_j=X_j(z)=X_{i,j}(z)=\sum_{m=0}^{\infty}p_{ij}^{(m)}z^m\]
 Expected number of arrivals in $j$ when starting in $i$:
 \[x_j=x_{i,j}=X_{j}(1)\]
  For MFB $M[j]$:
 \[Y_j=Y_j(z)=Y_{i,j}(z)=\sum_{m=0}^{\infty}p_{i,M[j]}^{(m)}z^m\]
 Expected number of arrivals in $M[j]$ when starting in $i$:
 \[y_j=y_{i,j}=Y_{j}(1)\]
On  edge $[i,j]$: 
\[\rho=\rho[i,j]=\frac{p[i,j]}{q[i,j]}; \quad \quad  n=n[i,j]\]

\begin{theorem} \label{t1}
 $y_k \ (k=0,1\dots,N)$ is the unique solution of $\sum_{j=0}^{N}u_{ij}y_j=-\delta(i,0)\ (i=0,1\dots N)$  where 
 \[u_{ij}=\left[ \frac{(1-\rho)\rho^n}{1-\rho ^{n+1}}\right]p_{[j,i]}^{*} \quad \quad (j<i,\ \rho\neq 1) \]
  \[u_{ij}=\left[ \frac{1-\rho}{1-\rho ^{n+1}}\right]p_{[j,i]}^{*} \quad \quad (j>i, \ \rho\neq 1)\]
  \[u_{ij}=\frac{p_{[j,i]}^{*}}{n+1} \quad \quad (j\neq i, \ \rho=1)\]
  \[u_{ii}=-1+\sum_{j<i, \rho\neq 1}\left[ \frac{\rho(1-\rho^n)}{1-\rho ^{n+1}}\right]p_{[i,j]}^{*} +\sum_{j>i, \rho\neq 1}\left[ \frac{1-\rho^n}{1-\rho ^{n+1}}\right]p_{[i,j]}^{*} 
  +\sum_{j\neq i, \rho= 1}\left[ \frac{n}{{n+1}}\right]p_{[i,j]}^{*}  \]
 \end{theorem}
 \begin{proof}
 CASE 1: $(z\neq 1) \lor (p\neq q)$. We prove the results of CASE 1 in 5 steps.\\
 Step 1:
 The random walk between $M[i]$ and $M[j] \ (0\leq i<j\leq N)$ is described by  considering the last step of the random walk:
  \begin{equation} \label{1}
 X_k=pzX_{k-1}+qzX_{k+1}
 \end{equation}
 where $X_k=X_{k}^{[i,j]}$, $p=p[i,j], q=q[i,j]$. Characteristic equation:
   \begin{equation} \label{2}
 qz\lambda^2-\lambda+pz=0
 \end{equation}
 with solutions $\lambda_1$ and $\lambda_2$ with $\lambda_1>\lambda_2$ and $\lambda_1\lambda_2=\frac{p}{q}=\rho$. So:
    \begin{equation}\label{3}
 X_k=a\lambda_1^{k}+b\lambda_2^k \quad (\lambda_1 > \lambda_2)
 \end{equation}
 Step 2: We express $a$ and $b$ in $Y_i$ and $Y_j$.
Focus on states 1 and $n=n[i,j]$ between $M[i]$ and $M[j] \ (0\leq i<j\leq N)$ and their neighbors. Considering the last step of the random walk we get:
 \begin{equation}\label{4}
 X_1=p_{[i,j]}^*zY_i+qzX_2
\end{equation}
\begin{equation}\label{5}
 X_n=pzX_{n-1}+p_{[j,i]}^*zY_j
 \end{equation}
 Using \eqref{3}, \eqref{4} and \eqref{5} we get:
 \begin{equation}\label{6}
(\lambda_2^{n+1}-\lambda_1^{n+1})a=\lambda_2^{n+1}\frac{p_{[i,j]}^*}{p}Y_i-\frac{p_{[j,i]}^*}{q}Y_j
 \end{equation}
 \begin{equation}\label{7}
(\lambda_2^{n+1}-\lambda_1^{n+1})b=\frac{p_{[j,i]}^*}{q}Y_j-\lambda_1^{n+1}\frac{p_{[i,j]}^*}{p}Y_i
 \end{equation}
 Step 3: We express $X_1$ and $X_n$ in $Y_i$ and $Y_j$. 
  Using \eqref{3} with $k=1$ and $k=n$ in combination with \eqref{6} and \eqref{7} gives:
 \begin{equation}\label{9}
q(\lambda_2^{n+1}-\lambda_1^{n+1})X_1=(\lambda_2-\lambda_1)p_{j,i}^*Y_j+(\lambda_2^{n}-\lambda_1^{n})p_{i,j}^*Y_i \quad (i<j)
 \end{equation}
 \begin{equation*}
p(\lambda_2^{n+1}-\lambda_1^{n+1})X_n=(\lambda_2-\lambda_1)(\lambda_1 \lambda_2)^n p_{i,j}^*Y_i+\lambda_1 \lambda_2(\lambda_2^{n}-\lambda_1^{n})p_{j,i}^*Y_j \quad (i<j)
 \end{equation*}
 We need the last formula for $j<i$, so we interchange $i$ and $j$:
 \begin{equation}\label{10}
p(\lambda_2^{n+1}-\lambda_1^{n+1})X_n=(\lambda_2-\lambda_1)(\lambda_1 \lambda_2)^n p_{j,i}^*Y_j+\lambda_1 \lambda_2(\lambda_2^{n}-\lambda_1^{n})p_{i,j}^*Y_i \quad (j<i)
 \end{equation}
 Step 4: Focus on $M[i]$ and its neighbors $X_1[i,j] \ (j>i)$ and $X_n[i,j] \ (j<i)$.
  Using \eqref{9} and \eqref{10} we get, considering the last step of the random walk:
 \begin{multline*}
 Y_i=\sum_{j>i}q[i,j]zX_1[i,j]+\sum_{j<i}p[i,j]zX_n[i,j]+\delta(i,0)=\\
 z\sum_{j>i}(\lambda_2^{n+1}-\lambda_1^{n+1})^{-1}\left[(\lambda_2-\lambda_1)p_{j,i}^*Y_j+(\lambda_2^{n}-\lambda_1^{n})p_{i,j}^*Y_i\right]+\\
 z\sum_{j<i}(\lambda_2^{n+1}-\lambda_1^{n+1})^{-1}\left[(\lambda_2-\lambda_1)(\lambda_1 \lambda_2)^n p_{j,i}^*Y_j+\lambda_1 \lambda_2(\lambda_2^{n}-\lambda_1^{n})p_{i,j}^*Y_i\right]+\delta(i,0)=\\
  zY_i \left\{ \sum_{j>i}\left [\frac{\lambda_2^n-\lambda_1^n}{\lambda_2^{n+1}-\lambda_1^{n+1}}\right ]p_{[i.j]}^*+\sum_{j<i}\left [\frac{\lambda_1\lambda_2(\lambda_2^n-\lambda_1^n)}{\lambda_2^{n+1}-\lambda_1^{n+1}}\right ]p_{[i,j]}^*\right\} +\\
 zY_j\left\{ \sum_{j>i}\left [\frac{\lambda_2-\lambda_1}{\lambda_2^{n+1}-\lambda_1^{n+1}}\right ]p_{[j,i]}^* +\sum_{j<i}\left [\frac{(\lambda_2-\lambda_1)(\lambda_1\lambda_2)^n }{\lambda_2^{n+1}-\lambda_1^{n+1}}\right ]p_{[j,i]}^*    \right\} + \delta(i,0)\\
\end{multline*}
Step 5:
If  $(z=1) \land (\rho>1)$ then  $\lambda_1=\rho; \lambda_2=1$. If $(z=1) \land (\rho<1)$ then $\lambda_1=1; \lambda_2=\rho$. We get the result by observing the coefficients of $y_i$ and $y_j$ and the constants. \\
CASE 2: $(z=1) \land (p=q)$ \\
  We can use  the same method as in CASE 1, but now with $x$ and $y$ instead of $X$ and $Y$, starting with $x_k=ak+b$, but we prefer a faster way by applying  l'Hospitals rule in the asymmetric case:
     \[\lim_{\rho \to 1} \frac{(1-\rho)\rho^n}{1-\rho ^{n+1}}=\frac{1}{n+1}=
  \lim_{\rho \to 1} \frac{1-\rho}{1-\rho ^{n+1}}\]
    \[\lim_{\rho \to 1} \frac{\rho(1-\rho^n)}{1-\rho ^{n+1}}=\frac{n}{n+1}=\lim_{\rho \to 1}\frac{1-\rho^n}{1-\rho ^{n+1}}\]
 \end{proof}

 \begin{theorem}\label{t2}
 Case $\rho \neq 1$:
  \begin{equation} 
  x_k=\frac{(1-\rho^k)\frac{p_{[j,i]}^*}{q}y_j+(\rho^k-\rho^{n+1})\frac{p_{[i,j]}^*}{p}y_i}{1-\rho^{n+1}}
  \end{equation}
  Case $\rho= 1$:
  \begin{equation} 
  x_k=\frac{kp_{[j,i]}^*y_j+(n+1-k)p_{[i,j]}^*y_i}{2(n+1)}
  \end{equation}
  \end{theorem}
  \begin{proof}
  Case $\rho \neq 1$:
 Use  \eqref{3}, \eqref{6} and \eqref{7} with $z=1$ and $\lambda_1=\rho; \lambda_2=1 (\rho>1)$ or $\lambda_1=1; \lambda_2=\rho (\rho<1)$. The result of $\rho=1$ is obtained by applying l'Hospitals rule for the asymmetric case.
  \end{proof}

 \begin{theorem}\label{t3}
 \begin{equation} \label{12}
 \sum_{j=0}^Np_{[j,j]}^*y_j=1
 \end{equation}
 \end{theorem}
 \begin{proof}
 Consider first $\rho \neq 1$. Using Theorem \ref{t1} we get: 
\begin{multline*}
\sum_{i=0}^Nu_{ij}=u_{jj}+\sum_{i<j}u_{ij}+\sum_{i>j}u_{ij}=-1+ \sum_{j<i}\left[ \frac{\rho(1-\rho^n)}{1-\rho ^{n+1}}\right]p_{[i,j]}^{*} + \\ \sum_{j>i}\left[ \frac{1-\rho^n}{1-\rho ^{n+1}}\right]p_{[i,j]}^{*}+
 \sum_{i<j}\left[ \frac{1-\rho}{1-\rho ^{n+1}}\right]p_{[j,i]}^{*}+\sum_{i>j} \left[ \frac{(1-\rho)\rho^n}{1-\rho ^{n+1}}\right]p_{[j,i]}^{*} 
\end{multline*}
Interchange $i$ and $j$ in the first two summations and add terms with $\sum_{i<j}$ and $\sum_{i>j}$: 
\begin{equation}\label{13}
\sum_{i=0}^Nu_{ij}=-1+\sum_{i\neq j}p_{[j,i]}^*=-p_{[j,j]}^*
 \end{equation}
 We also have (Theorem \ref{t1}) $\sum_{j=0}^{N}u_{ij}y_j=-\delta(i,0)\ (i=0,1\dots N)$, so $\sum_{i=0}^N\sum_{j=0}^{N}u_{ij}y_j=-1.$ 
 Using \eqref{13} we get:  $\sum_{i=0}^N\sum_{j=0}^{N}u_{ij}y_j=\sum_{j=0}^N\left[\sum_{i=0}^{N}u_{ij}\right ]y_j=-\sum_{j=0}^Np_{[j,j]}^*y_j.$
The symmetric case $\rho=1$ proceeds along the same lines.
\end{proof}
\subsection{Expected time until absorption}
Let $t_k$ be the expected time until absorption when starting in $M[k] \ (k=0,1,\dots,N)$.

\begin{theorem} \label{t4}
 $t_k \ (k=0,1\dots,N)$ is the unique solution of $\sum_{j=0}^{N}v_{ij}t_j=\tau_i \ (i=0,1\dots N)$  where 
  \[v_{ij}=\left[ \frac{(1-\rho)\rho^n}{1-\rho ^{n+1}}\right]p_{[i,j]}^{*} \quad \quad (j>i,\ \rho\neq 1) \]
  \[v_{ij}=\left[ \frac{1-\rho}{1-\rho ^{n+1}}\right]p_{[i,j]}^{*} \quad \quad (j<i, \ \rho\neq 1)\]
  \[v_{ij}=\frac{p_{[i,j]}^{*}}{n+1} \quad \quad (j\neq i, \ \rho=1)\]
  \[v_{ii}=-1+\sum_{j<i, \rho\neq 1}\left[ \frac{\rho(1-\rho^n)}{1-\rho ^{n+1}}\right]p_{[i,j]}^{*} +\sum_{j>i, \rho\neq 1}\left[ \frac{1-\rho^n}{1-\rho ^{n+1}}\right]p_{[i,j]}^{*} 
  +\sum_{j\neq i, \rho= 1}\left[ \frac{n}{{n+1}}\right]p_{[i,j]}^{*}  \]
  \begin{multline*}
  \tau_i=-1+ \sum_{j<i, \rho\neq 1}\left[ \frac{n-(n+1)\rho+\rho^{n+1}}{(p-q)(1-\rho ^{n+1})}\right]p_{[i,j]}^{*} +\\
  \sum_{j>i, \rho\neq 1}\left[ \frac{1-(n+1)\rho^n+n\rho ^{n+1}}{(p-q)(1-\rho^{n+1})}\right]p_{[i,j]}^{*}-\sum_{j\neq i, \rho= 1}np_{[i,j]}^{*}
  \end{multline*}  
 \end{theorem}
\begin{proof}
Step 1: Let $m_k=m_k[i,j] \ (k=1,2,\dots,n[i,j])$ be the expected time until absorption when starting on edge $[i,j]$ in state $k \ (k=1,2,\dots,n[i,j])$. We have, considering the next step in the random walk:
\begin{equation*}
m_k=p(m_{k+1}+1)+q(m_{k-1}+1)=pm_{k+1}+qm_{k-1}+1 \quad (k=1,2,\dots,n[i,j])
\end{equation*}
with general solution (case $\rho\neq 1$):
\begin{equation}\label{15}
m_k=a\rho^{-k}+b-\frac{k}{p-q} \quad (k=0,1,\dots,n[i,j]+1)
\end{equation}
Step 2: 
We express $a$ and $b$ in $t_i$ and $t_j$ using \eqref{15} with $k=0$ and $k=n+1$: $m_0=a+b=t_i$ and $m_{n+1}=a\rho^{-n-1}+b-\frac{n+1}{p-q}=t_j$ gives $(1-\rho^{n+1})a=-\rho^{n+1}t_i+\rho^{n+1}t_j+\rho^{n+1}\left(\frac{n+1}{p-q}\right )$ and $(1-\rho^{n+1})b=t_i-\rho^{n+1}t_j-\rho^{n+1}\left(\frac{n+1}{p-q}\right )$ \\
Step 3: 
Using the expressions for $a$ and $b$ we get, using \eqref{15}:
\begin{equation*}
m_1=\frac{(1-\rho^n)t_i+(\rho^n-\rho^{n+1})t_j}{1-\rho^{n+1}}+\frac{(n+1)(\rho^n-\rho^{n+1})}{(p-q)(1-\rho^{n+1})}-\frac{1}{p-q} \quad (i< j)
\end{equation*}

\begin{equation*}
m_n=\frac{(1-\rho)t_i+(\rho-\rho^{n+1})t_j}{1-\rho^{n+1}}+\frac{(n+1)(\rho-\rho^{n+1})}{(p-q)(1-\rho^{n+1)}}-\frac{n}{p-q} \quad (i< j)
\end{equation*}
We need the last formula with $i>j$, so by interchanging $i$ and $j$ we get:

\begin{equation*}
m_n=\frac{(1-\rho)t_j+(\rho-\rho^{n+1})t_i}{1-\rho^{n+1}}+\frac{(n+1)(\rho-\rho^{n+1})}{(p-q)(1-\rho^{n+1)}}-\frac{n}{p-q} \quad (i> j)
\end{equation*}
Step 4. Consider $\rho=1$ on edge $[i,j]$: $m_k=\frac{1}{2}m_{k-1}+\frac{1}{2}m_{k+1}+1$ gives
\begin{equation*}
m_k=ak+b-k^2 \quad (k=0,1,\dots,n[i,j]+1)
\end{equation*}
Along the same lines as in case $\rho\neq 1$ we get:
\begin{equation*}
m_1=\frac{nt_i+t_j}{n+1}+n
 \quad (i< j)
\end{equation*}
\begin{equation*}
m_n=\frac{nt_i+t_j}{n+1}+n
 \quad (i> j)
\end{equation*}
The same formula are found by applying l'Hospitals rule twice in case of $\rho\neq 1$.
Step 5: Substituting the values of $m_1$ and $m_n$ and considering the next step of the random walk:
\begin{multline*}
t_i=\sum_{j>i}p_{[i,j]}^*(m_1[i,j]+1)+\sum_{j<i}p_{[i,j]}^*(m_n[i,j]+1)+p_{[i,i]}^*.1=\\
1+\sum_{j>i}p_{[i,j]}^*m_1[i,j]+\sum_{j<i}p_{[i,j]}^*m_n[i,j]=\\
1+\sum_{j>i,\rho \neq 1}p_{[i,j]}^*\left [\frac{(1-\rho^n)t_i+(\rho^n-\rho^{n+1})t_j}{1-\rho^{n+1}}+\frac{(n+1)(\rho^n-\rho^{n+1})}{(p-q)(1-\rho^{n+1})}-\frac{1}{p-q}\right]+\\
\sum_{j<i,\rho \neq 1}p_{[i,j]}^*\left[\frac{(1-\rho)t_j+(\rho-\rho^{n+1})t_i}{1-\rho^{n+1}}+\frac{(n+1)(\rho-\rho^{n+1})}{(p-q)(1-\rho^{n+1)}}-\frac{n}{p-q}\right]+\\
\sum_{j>i,\rho = 1}p_{[i,j]}^*\left[\frac{nt_i+t_j}{n+1}+n\right ]+\sum_{j<i,\rho = 1}p_{[i,j]}^*\left [\frac{nt_i+t_j}{n+1}+n\right ]               
\end{multline*}
where $ i=0,1,\dots,N$. By observing the coefficients of $t_i$ and $t_j$ and the constants we obtain the result.
\end{proof}
\section{Examples of multiple function barrier graphs}
\subsection{A finite star graph}
We consider a star graph where the starting state is in the center $M[0]$. In a finite star graph we demand:
\[\sum_{i=0}^Np_{[0,i]}^*=1 \quad \quad \prod_{i=1}^Np_{[0,i]}^*\neq 0\]
 \[p_{[i,0]}^*+p_ {[i,i]}^*=1 \quad (i=1,2,\dots,N)\] 
Example with $N=4$:
 \graphicspath{file:///C:/Users/Theo/Desktop/}
\begin{center}
\includegraphics[scale=0.4]{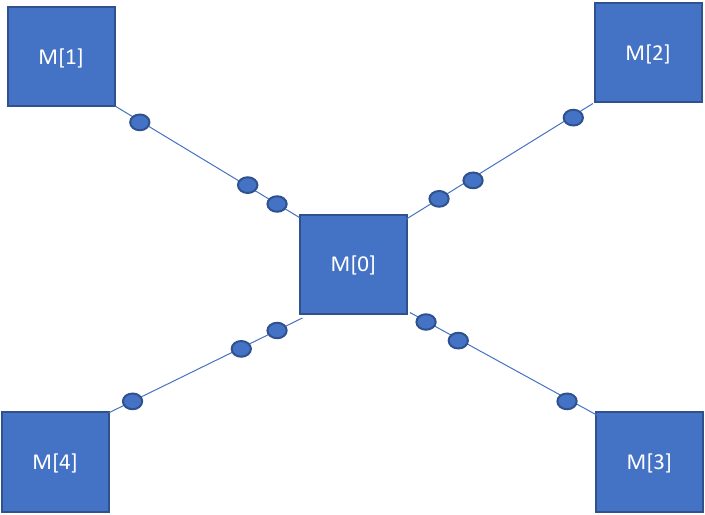} 
\end{center}
 We use the notation: $p_i=p[0,i], \quad q_i=q[0,i], \quad \rho_i=\frac{p_i}{q_i} \quad (i=1,2,\dots,N)$
\subsubsection{Expected number of arrivals}
To obtain the expected number of arrivals in the finite star graph, we use Theorem \ref{t1} : $\sum_{j=0}^{N}u_{ij}y_j=-\delta(i,0)\ (i=0,1\dots N)$. 
We get (case $\rho\neq 1$): 
\begin{equation}\label{16}
y_i=\left[-\frac{u_{i0}}{u_{ii}}\right]y_0=\left[\frac{(1-\rho_i)\rho_i^np_{[0,i]}^*}{(1-\rho_i)+\rho_i(1-\rho_i^n)p_{[i,i]}^*}\right]y_0 \quad (i=1,2,\dots,N)
\end{equation}
When $\rho=1$ we get:
\begin{equation} \label{17}
y_i=\frac{p_{[0,i]}^*}{1+n_ip_{[i,i]}^*}y_0 \quad (i=1,2,\dots,N)
\end{equation}
Instead of the first equation of $\sum_{j=0}^{N}u_{ij}y_j=-\delta(i,0)\ (i=0,1\dots N)$ we use the result of Theorem \ref{t3}: $\sum_{j=0}^Np_{[j,j]}^*y_j=1$, which leads to:
\begin{equation}\label{18}
y_0=\frac{1}{p_{[0,0]}^*+\sum_{i=1, \rho_i\neq 1}^N\left[\frac{(1-\rho_i)\rho_i^np_{[0,i]}^*}{(1-\rho_i)+\rho_i(1-\rho_i^n)p_{[i,i]}^*}\right]+\sum_{i=1, \rho_i=1}^N\left[\frac{p_{[0,i]}^*}{1+n_ip_{[i,i]}^*}\right]} 
\end{equation}

\subsubsection{Mean absorption time with absorbing barriers}
We consider a star graph with starting point $M[0]$ in the center and all other MFB's are absorbing:
\[ \sum_{i=0}^N p_{[0,i]}^*=1;\quad \prod_{i=1}^Np_{[0,i]}^*\neq 0 ;\quad p_{[i,i]}^*=1 \ (i=1,2,\dots,N)   \]
We use notation:
\[p_i=p[0,i], \ q_i=q[0,i], \ n_i=n[0,i], \ \rho_i=\frac{p_i}{q_i} \ (i=1,2,\dots,N)\]
Let for $i=1,2,\dots,N$:
\[ \alpha_i=\frac{1-\rho_i^{n_i}}{1-\rho_i^{n_i+1}} ; \quad \quad \beta_i=\frac{1-(1+n_i) \rho_i^{n_i}+n_i \rho_i^{n_i+1}}{(p_i-q_i)(1- \rho_i^{n_i+1})}\quad \quad (\rho_i\neq 1)   \]
\[ \alpha_i=\frac{n_i}{n_i+1}; \quad \quad \beta_i=-n_i \quad \quad (\rho_i=1)\]
Theorem \ref{t4} gives:
\[t_0=\frac{\tau_0}{v_{00}}=\frac{1-\sum_{i=1}^{N}\beta_ip_{[0,i]}^*}{1-\sum_{i=1}^N \alpha_ip_{[0,i]}^*}\]
\subsection{An infinite star graph with absorbing barriers}
In this subsection we consider an infinite star graph where all barriers (except the start position $M[0]$) are absorbing: $p_{[i,i]}^*=1 \quad (i=1,2,\dots)$.
\begin{theorem}
\begin{equation}\label{19}
y_0=\frac{1}{1-\sum_{i=1, \rho_i\neq 1}^\infty \left[\frac{1-\rho_i^n}{1-\rho_i^{n+1}} \right]p_{[0,i]}^*-\sum_{i=1,\rho_i=1}^{\infty}\left[\frac{n}{1+n}\right]p_{[0,i]}^*} \quad (n=n_i)
\end{equation}
\begin{equation}
y_i=\frac{(1-\rho_i)\rho_i^np_{[0,i]}^*}{1-\rho_i^{n+1}} y_0 \ \  (\rho_i\neq 1;\ n=n_i; \ i=1,2,\dots)
\end{equation}
\begin{equation}
y_i=\frac{p_{[0,i]}^*}{1+n}y_0 \ \ (\rho_i= 1; \ n=n_i; \ i=1,2,\dots)\end{equation}
\end{theorem}
\begin{proof}
Use \eqref{16},\eqref{17} and \eqref{18} with $p_{[i,i]}^*=1 (i=1,2,\dots,N)$ and rewrite \eqref{18} to the form in \eqref{19} by using $p_{[i,i]}^*=1.$ Finally note that
$\sum p_{[0,i]}^*$ is a majorant of both  $\sum \left[\frac{n}{1+n}\right]p_{[0,i]}^*$ and $\sum \left[\frac{1-\rho_i^n}{1-\rho_i^{n+1}} \right]p_{[0,i]}^*$
\end{proof}

\subsection{A positive oriented finite cycle graph}
We have $N+1$ barriers in the finite cycle graph: $M[0],M[1],\dots,M[N]$.  We start in $M[0]$. When the random walk is in $M[i]$ then absorption can happen or we move one step in the direction of $M[i+1]$. The state space is mod $(N+1)$: when arriving in $M[N]$, there can be a step in the direction of $M[N+1]$ where $M[N+1]=M[0]$. We have:
\begin{equation} \label{23}
p_{[i,i]}^*+p_{[i,i+1]}^*=1 \quad (i=0,1,\dots,N)
\end{equation}
Example with $N=3$:
\graphicspath{file:///C:/Users/Theo/Desktop/}
\begin{center}
\includegraphics[scale=0.4]{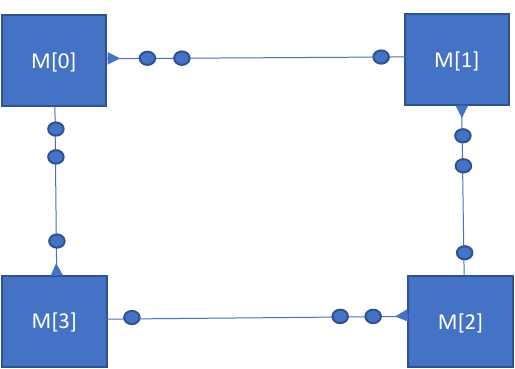} 
\end{center}
We use the notation: $\rho_i=\rho[i,i+1]; \quad n_i=n[i,i+1]$.
\subsubsection{Expected number of arrivals}
\begin{theorem}
\begin{equation}\label{24}
y_k=\frac{\prod_{i=1}^{k}M_i}{\sum_{m=0}^{N} p_{[m,m]}^*\prod_{i=1}^{m}M_i}
 \quad (k=0,1,\dots,N) \quad (\prod_{i=1}^{0}M_i=1)
\end{equation}
case $\rho \neq 1$:
\begin{equation}\label{25}
M_i=\frac{p_{[i-1,i]}^*}{1+\left[\frac{1-\rho^{n}}{\rho^{n}(1-\rho)}   \right]p_{[i,i]}^*}
 \quad (i=1,2,\dots,N)\quad \rho=\rho_i; \quad n=n_i
\end{equation}
case $\rho =1$:
\begin{equation}\label{26}
M_i=\frac{p_{[i-1,i]}^*}{1+n_ip_{[i,i]}^*}
 \quad (i=1,2,\dots,N)
\end{equation}
\end{theorem}
\begin{proof}
We use \eqref{23} and Theorem \ref{t1} to obtain:
$y_i=-\frac{u_{i,i-1}}{u_{ii}}y_{i-1}
\quad (i=1,2,\dots,N)$.
Let $M_i=-\frac{u_{i,i-1}}{u_{ii}}$, $\ \rho=\rho_i$ and $n=n_i$ then $M[i]=\frac{p_{[i-1,i]}^*}{1 +\left [\frac{1-\rho^n}{(1-\rho)\rho^n} \right ]  p_{[i,i]}^*} \ (\rho \neq 1)$ and $M[i]=\frac{p_{[i-1,i]}^*}{1+np_{[i,i]}^*} \ (\rho=1)$. Because of $y_i=M_iy_{i-1} \ (i=1,2,\dots,N)$ we have $y_k=\left( \prod_{i=1}^kM_i \right )y_0 \ \ (k=0,1,\dots,N)$ where we define $\prod_{i=1}^0M_i=1$. By theorem \ref{t3} we have $\sum_{j=0}^Np_{[j,j]}^*y_j=1$. Combining the last two results gives \eqref{24}.
\end{proof}

\subsubsection{Mean absorption time}

We use notation:
\[p_i=p[i,i+1], \ q_i=q[i,i+1], \ n_i=n[i,i+1], \ \rho_i=\frac{p_i}{q_i} \ (i=0,1,\dots,N)\]
Let for $i=0,1,\dots,N$:
\[ \alpha_i=\frac{1-\rho_i^{n_i}}{1-\rho_i^{n_i+1}} ; \quad \quad \beta_i=\frac{1-(1+n_i) \rho_i^{n_i}+n_i \rho_i^{n_i+1}}{(p_i-q_i)(1- \rho_i^{n_i+1})};\quad \quad \gamma_i=\frac{(1-\rho_i)\rho_i^{n_i}}{1-\rho_i^{n_i+1}} \quad \quad(\rho_i\neq 1)   \]
\[ \alpha_i=\frac{n_i}{n_i+1}; \quad \quad \beta_i=-n_i; \quad \quad \gamma_i=\frac{1}{n_i+1} \quad \quad (\rho_i=1)\]

Using Theorem \ref{t4} we obtain:
\[v_{i,i+1}=\gamma_{i}p_{[i,i+1]}^*; \quad    v_{i,i}=-1+\alpha_i p_{[i,i+1]}^*; \quad \tau_{i}=-1+\beta_ip_{[i,i+1]}^* \] 
\[t_{i+1}=\frac{\tau_i-v_{i,i}t_i}{v_{i,i+1}}=\lambda_it_i+\mu_i; \quad \lambda_i=\frac{-v_{i,i}}{v_{i,i+1}};\quad \mu_i=\frac{\tau_i}{v_{i,i+1}} \quad (i=0,1,\dots,N)   \]
\[t_{k+1}=t_0 \prod_{i=0}^{k}\lambda_i+\sum_{i=0}^{k-1}\mu_i \prod_{j=i+1}^{k} \lambda_j \quad (k=0,1,\dots,N)\]
\[t_0=t_{N+1}=\frac{ \sum_{i=0}^{N-1}\mu_i \prod_{j=i+1}^{N} \lambda_j  }{1-\prod_{i=0}^{N}\lambda_i}   \]
\subsection{A positive oriented infinite cycle graph}
\begin{theorem}
CASE $\rho \neq 1$: 
If $\sum_{i=0}^{\infty} p_{[i,i]}^*$ converges and  \\ $\liminf_{i\to\infty}\rho_i>0, \limsup_{i\to\infty, \rho_i< 1}\rho_i<1, \liminf_{i\to\infty, \rho_i> 1}\rho_i>1 $ , then:
\begin{equation}\label{27}
y_k=\frac{\prod_{i=1}^{k}M_i}{\sum_{m=0}^{\infty} p_{[m,m]}^*\prod_{i=1}^{m}M_i}
 \quad (k=0,1,\dots) \quad (\prod_{i=1}^{0}M_i=1)
\end{equation}
where $M_i$ is given by $\eqref{25}$. \\
CASE $\rho=1$: We get \eqref{27} when $\sum_{i=0}^{\infty}n_i p_{[i,i]}^*$ converges, where $M_i$ is given by $\eqref{26}$.
\end{theorem}
\begin{proof}
We use: If  $0\leq \omega_{i}<1 $ then $ \prod _{i=1}^{\infty }(1-\omega_{i})$ converges to a non-zero number if and only if $\sum _{i=1}^{\infty }\omega_{i}$ converges. We define a relation $\backsim$ between two sequences $\sum a_i$ and $\sum b_i$: $\sum a_i \backsim \sum b_i$ if and only if both sequences converges.
Let $\omega_i=1-M_i$. \\
CASE $\rho \neq 1$:
Using \ref{25} we get: 
\begin{equation*}
0<\omega_i=\frac{(1-\rho_i^n)p_{[i,i]}^*+(1-\rho_i)\rho_i^np_{[i-1,i-1]}^*}{(1-\rho_i^n)p_{[i,i]}^*+(1-\rho_i)\rho_i^n}<1
\end{equation*}
If only a finite number of the $\rho_i$ are in the neighbor of 0 and 1 \\ ( $\liminf_{i\to\infty}\rho_i>0 , \limsup_{i\to\infty, \rho_i< 1}\rho_i<1, \liminf_{i\to\infty, \rho_i> 1}\rho_i>1)$ then by the comparison criterium and $\lim_{i \to \infty}p_{[i,i]}^*=0$: 
\begin{multline*}
\sum \omega_i=\sum \frac{(1-\rho_i^n)p_{[i,i]}^*+(1-\rho_i)\rho_i^np_{[i-1,i-1]}^*}{(1-\rho_i^n)p_{[i,i]}^*+(1-\rho_i)\rho_i^n}   \backsim \\ 
\sum \frac{(1-\rho_i^n)p_{[i,i]}^*+(1-\rho_i)\rho_i^np_{[i-1,i-1]}^*}{(1-\rho_i)\rho_i^n}\backsim \sum p_{[i,i]}^* 
\end{multline*}
 CASE $\rho=1$: Using \ref{26} we get : 
$0<\omega_i=\frac{n_ip_{[i,i]}^*+p_{[i-1,i-1]}^*}{n_ip_{[i,i]}^*+1}<1$
By the comparison criterium and $\lim_{i \to \infty}n_i p_{[i,i]}^*=0$: 
$\sum  \omega_i=\sum \frac{n_ip_{[i,i]}^*+p_{[i-1,i-1]}^*}{n_ip_{[i,i]}^*+1}\backsim \sum [n_ip_{[i,i]}^*+p_{[i-1,i-1]}^*]\backsim \sum n_ip_{[i,i]}^*$ (because of $n_i\geq 1$).
\end{proof}
%
\section{References}
  author       = {Weesakul, B.},
  date         = {1961},
  journaltitle = {Ann. Math. Statist., 32 765-769},
  title        = {The random walk between a reflecting and an absorbing barrier} \\ \\
  author = {Feller, W.},
  date   = {1968},
  title  = {An Introduction to probability theory and its applications, (third edition) Vol. 1 John Wiley, New York} \\ \\
  author       = {Lehner, G.},
  date         = {1963},
  journaltitle = {Ann. Math.Stat., 34, 405-412},
  title        = {One-dimensional random walk with a partially reflecting barrier}\\ \\
  author       = {Gupta, H. C.},
  date         = {1966},
  journaltitle = {Journ. Math Sci,1 18-29},
  title        = {Random walk in the presence of a multiple function barrier}\\ \\
  author       = {Dua, S. and Khadilkar, S. and Sen, K.},
  date         = {1976},
  journaltitle = {J. Appl. Prob., 13 169-175},
  title        = {A modified random walk in the presence of partially reflecting barriers}\\ \\
  author       = {Percus, O. E.},
  date         = {1985},
  journaltitle = {Adv. Appl. Prob. 17 594-606},
  title        = {Phase transition in one-dimensional random walk with partially reflecting boundaries}\\ \\
  author       = {El-Shehawey, M. A.},
  date         = {2000},
  journaltitle = {J. Phys. A: Math. Gen., 33, 9005-9013.},
  title        = {Absorption probabilities for a random walk between two partially absorbing boundaries. I}\\ \\
  author       = {Burioni, R. and Cassi, D.},
  date         = {2005},
  journaltitle = {J. Phys.A: Math. Gen., 38, R45-R78.},
  title        = {Random walks on graphs: ideas, techniques and results} \\ \\
  author       = {Durhuus, B. and Jonsson, T. and Wheater, J.},
  date         = {2006},
  journaltitle = {J. Phys. A: Math. Gen.,39, 1009-1038},
  title        = {Random walks on combs} \\ \\
\end{document}